\documentclass[11pt]{article}

\usepackage{CJK,indentfirst,amscd,amsmath,amsfonts,amssymb,amsthm,cite}
\usepackage{graphicx,epstopdf,subfigure,tikz,xcolor,verbatim}
\usepackage[all]{xy}
\usepackage[perpage,symbol*]{footmisc}

\newtheorem{theorem}{Theorem}
\newtheorem{lemma}{Lemma}

\newtheorem{corollary}{Corollary}
\newtheorem{conjecture}{Conjecture}

\newcommand{\eproof}{\hspace*{\fill}{$\blacksquare$}}

\begin{document}
\pagestyle{myheadings}

\title{\bf On the algebraic bipartite graphs $D(4s,q)$ \footnote{This work was supported by the Natural Science Foundation of China (No. 61977056)
and Natural Science Foundation of the Jiangsu Higher Education Institutions of China (25KJB110013).}}
	\author{ Chen Zhang, Xiaoyan Cheng, Yuansheng Tang\footnote{Corresponding author.} \\
{\it\small School of Mathematical Sciences, Yangzhou University, Jiangsu, China}\\
and Ming Xu
{\it\small \footnote{Email addresses: cz12024@163.com(C. Zhang), xycheng@yzu.edu.cn(X. Cheng), ystang@yzu.edu.cn(Y. Tang), mxu@szcu.edu.cn(M.Xu).}}\\
{\it\small Suzhou City University,Jiangsu, China}
           }
	\date{}
	\maketitle

\noindent\small{\bf Abstract:}
Lazebnik-Ustimenko graphs $D(k,q)$ constitute a notable family of $q$-regular, edge-transitive bipartite graphs characterized by high symmetry and large girth, with wide applications in algebraic graph theory and coding theory. The precise girth of $D(k,q)$ is governed by a celebrated conjecture proposed by F\"{u}redi et al., which remains unresolved in general cases. Although prior works have validated this conjecture for several specific parameter regimes and established fundamental structural properties of graph girth, a systematic characterization of girth cycles for even-indexed subfamilies of \(D(k,q)\) is still absent. This paper focuses on the graph subclass $D(4s,q)$ and conducts structural analysis via its isomorphic bipartite model $\Lambda_{k,q}$. We derive an explicit necessary and sufficient condition for the existence of length-$(4s+4)$ backtrackless circuits in $\Lambda_{4s,q}$. Using this criterion, we completely characterize all girth cycles of $\Lambda_{8,q}$ and verify the F\"{u}redi-Lazebnik-Ustimenko girth conjecture for the parameter family $k=2(p^r+1)p^t-4$ over finite fields of odd characteristic $p\geq 3$ with $q=p^m>3$. The results further enrich the structural theory of symmetric large-girth bipartite graphs and extend the valid parameter range of this classic girth conjecture.
\vspace{1ex}


\section{Introduction}

Throughout this paper, all graphs considered are finite, undirected, and simple, containing no loops or multiple edges. A graph $G$ is edge-transitive if for any two edges in $E(G)$, there exists an automorphism of $G$ mapping the endpoints of one edge to those of the other. A backtrackless walk of length $n$ in $G$ is a vertex sequence where every pair of consecutive vertices is adjacent and no immediate backward step occurs. A closed backtrackless walk is referred to as a backtrackless circuit. For any non-forest graph $G$, the girth $g(G)$ is defined as the length of its shortest backtrackless circuit. Accordingly, every girth-attaining circuit is a simple cycle consisting of exactly $g(G)$ distinct vertices.

Graphs with large girth and high symmetry are fundamental research objects in extremal graph theory, finite geometry, coding theory, cryptography, network design, and quantum computation (cf. \cite{Wenger91}-\cite{Dehghan20}).
Let $q$ be a prime power of characteristic $p$. For $k\geq 2$, the $q$-regular edge-transitive bipartite graph $D(k,q)$, constructed by Lazebnik and Ustimenko \cite{Lazebnik93} and widely studied in subsequent literature (cf. \cite{Lazebnik93}-\cite{XCT2023}), admits a rigorous girth lower bound formulated as
\begin{align*}
g_k=\left\{
\begin{array}{ll}
  k+4, & \text{if }k\text{ is even}, \\
  k+5, & \text{if }k\text{ is odd}.
\end{array}
\right.
\end{align*}
The exact girth value of $D(k,q)$ was conjectured by F\"{u}rediet al. \cite{Furedi95}, as formally stated below.
   \begin{conjecture}\label{conmain}
     For all integer $k\geq 2$ and prime power $q\geq4$, the girth of $D(k,q)$ satisfies $g(D(k,q))=g_k$.
   \end{conjecture}

A notable reduction property facilitates the verification of Conjecture 1: the validity of the conjecture for even parameters $k=2t$ implies its validity for the adjacent odd parameters $k=2t-1$, since $g(D(2t,q))=2t+4$ forces $g(D(2t-1,q))=2t+4$. In \cite{Furedi95}, Conjecture~\ref{conmain} was verified for all pairs $(k,q)$ such that $g_k/2$ divides $q-1$. The proof relies on constructing an automorphism $\sigma$ of $D(k,q)$ of order $g_k/2$, for which some edge $e$ and its image $\sigma(e)$ jointly form a backtrackless walk of length four (cf. \cite{LazebnikSunWang17}).
Further progress was achieved in \cite{Cheng14} and \cite{Cheng16}, where explicit closed-form expressions for backtrackless walks were exploited to confirm Conjecture~\ref{conmain} in two important cases:
when $g_k/2$ is a pure power of the characteristic $p$ and when $g_k/2$ factors as a product of a divisor of $q-1$ and a $p$-power, respectively.

The original automorphism-based approach in \cite{Furedi95} was significantly refined by Sun \cite{Sun17}. By constructing automorphisms of order $g_k/4$ and utilizing length-4 backtrackless walks that generate length-8 orbits under automorphic actions, Sun confirmed the conjecture for odd primes $p$ and parameters $k=4hp^t-5$, where $h$ divides $q-1$. This result was independently reproved via an alternative method in \cite{XCT2023}, which additionally established three new structural relations for all integer $t\geq 1$:
(a) $g(D(4t+2,q))=g(D(4t+1,q))$;
(b) $g(D(4t+3,q))=g_{4t+3}$ whenever $g(D(2t,q))=g_{2t}$;
(c) $g(D(8t,q))=g_{8t}$ whenever $g(D(4t-2,q))=g_{4t-2}$.

Additionally, a backtrackless circuit lifting technique proposed in \cite{Sun17} generates infinite valid parameter families for Conjecture \ref{conmain}. Specifically, if the conjecture holds for $k=n$, it extends to the infinite family $\{g_np^s-5:s\geq 1\}$ when $n\equiv 3\pmod{4}$, and to $\{g_np^s-4:s\geq 1\}$ for all remaining cases (cf. \cite{CZTX26}).

In addition to their theoretical value, bipartite graphs with large girth find extensive applications in scientific and engineering modeling. Their two disjoint vertex partitions naturally represent system variables and local constraint rules acting on neighboring units, which renders such graphs particularly useful in coding theory. A typical application lies in the Tanner graph representation of low-density parity-check (LDPC) codes, where variable nodes correspond to code symbols and check nodes correspond to parity-check constraints (cf. \cite{Kim04,Yan11,Polak13}). Accordingly, explicit characterization of girth cycles in large-girth bipartite graphs is of both theoretical significance and practical relevance. In recent work \cite{XCT23}, complete classifications of girth cycles for $D(k,q)$ were obtained for low parameter regimes: $3\leq k\leq 5$ for $q>3$ and $3\leq k\leq 14$ for $q=3$.

This paper focuses on the graph family $D(4s,q)$. Following the standard framework in \cite{Cheng14,Cheng16,XCT2023}, we adopt the isomorphic bipartite model $\Lambda_{k,q}$ (rather than the original definition of $D(k,q)$) as the research object. Defined over $\mathbb{F}_q^{k+1}$, \(\Lambda_{k,q}\) has vertex sets specified as follows.
The left vertex set is
\[
L_k=\{[l]=(l_0,l_1,\ldots,l_k)\in \mathbb{F}_q^{k+1}: l_1=l_2\},
\]
and the right vertex set is
\[
R_k=\{\langle r\rangle=(r_0,r_1,\ldots,r_k)\in \mathbb{F}_q^{k+1}: r_1=0\}.
\]
A left vertex $[l]\in L_k$ is adjacent to a right vertex $\langle r\rangle\in R_k$ in $\Lambda_{k,q}$ if and only if for all $2\le i\le k$,
\begin{equation}\label{eq:adjacency}
l_i+r_i=
\begin{cases}
r_0l_{i-2}, & \text{if } i\equiv 2,3(\bmod 4),\\
l_0r_{i-2}, & \text{if } i\equiv 0,1(\bmod 4).
\end{cases}
\end{equation}
The isomorphism $\Lambda_{k,q}\cong D(k,q)$ is well-established in \cite{Cheng14}, which implies that $\Lambda_{k,q}$ inherits the edge-transitivity of $D(k,q)$ and satisfies the uniform girth lower bound $g(\Lambda_{k,q})\ge k+4$.
The contributions of this paper are twofold. First, we derive an explicit necessary and sufficient condition for the existence of length-$(4s+4)$ backtrackless circuits in $\Lambda_{4s,q}$. Second, based on this criterion, we fully characterize all girth cycles of $\Lambda_{8,q}$ and verify Conjecture \ref{conmain} for the parameter family $k=2(p^r+1)p^t-4$ under the conditions $p>2$ and $q>3$.

The remainder of this paper is organized as follows. Section \ref{sec:backtrackless-walks} reviews relevant homogeneous polynomials and derives simplified coordinate formulas for backtrackless walks in $\Lambda_{k,q}$. Section \ref{sec:backtrackless-Circuits} establishes the existence criterion for length-$(4s+4)$ backtrackless circuits in $\Lambda_{4s,q}$. Section \ref{sec:xtty} constructs explicit girth cycles, completes the full characterization of girth cycles of $\Lambda_{8,q}$, and classifies all $\Lambda_{4s,q}$ graphs that admit girth cycles with the symmetric pattern $(x,1,\dots,1,y)$.

\section{Backtrackless Circuits in $\Lambda_{k,q}$}\label{sec:backtrackless-walks}
At first, we introduce a class of homogeneous polynomials in several indeterminates which were defined in \cite{Cheng14}. For any tuple $\alpha=(\omega_1,\ldots,\omega_m)$ over $\mathbb{F}_q^*$, let
$$\rho_0(\alpha)=\omega_1\cdots\omega_m$$
and, for $1\leq s\leq \lfloor \frac{m}{2} \rfloor$, let $\rho_s(\alpha)$ denote the homogeneous polynomial of order $m-2s$ defined by
$$\rho_s(\alpha)=\sum_{1\leq i_1 <\cdots < i_s \leq m-s}\frac{\Pi_{j=1}^m\omega_j}{\Pi_{j=1}^s\omega_{i_j+j-1}\omega_{i_j+j}},$$
where each term in the summation is a product of the remaining entries of the tuple $\alpha$ under a deletion of $s$ disjoint pairs $(\omega_i,\omega_{i+1}).$ If $m < 2s$ or $s<0,$ $\rho_s(\alpha)$ is defined as $0$. For the null sequence $\eta$, we define $\rho_s(\eta)$ as
\[\rho_s(\eta)=
    \begin{cases}
    1 \quad \text{if }s=0,\\
    0 \quad \text{if }s\neq0.
    \end{cases}
\]
The following three lemmas can be derived easily from the definition of the homogeneous polynomial $\rho_s(\alpha)$ (cf. \cite{Cheng14,Cheng16,XCT23}).

\begin{lemma}\label{lem:rho}
For any tuples $\alpha,\beta$ over $\mathbb{F}_q^*$ and $\omega_1,\omega_2 \in \mathbb{F}_q^*,$
\begin{align}\label{eq:rho-k-1-2}
&\rho_k(\alpha,\omega_1,\omega_2,\beta)=\sum_i\rho_i(\alpha,\omega_1)
\rho_{k-i}(\omega_2,\beta)+\sum_j\rho_j(\alpha)\rho_
{k-1-j}(\beta).
\end{align}
In particular, we have
\begin{align}\label{eq:rho1}
\rho_s(\omega_1,\omega_2,\beta)=\rho_{s-1}(\beta)+\omega_1\rho_s(\omega_2,\beta),\\
\label{eq:rho2}
\rho_s(\alpha,\omega_1,\omega_2)=\rho_{s-1}(\alpha)+\omega_2\rho_s(\alpha,\omega_1).
\end{align}
\end{lemma}
\begin{lemma}\label{lemma20}
Let $\alpha=(\omega_1,\ldots,\omega_{m})$ be an $m$-tuple over $\mathbb{F}_q^*$.
\begin{enumerate}
  \item If $m=2n$ is even, then for  $0\leq j\leq n$ we have
\begin{align}\label{eq:rho-sum}
&\rho_{n-j}(\alpha)=\sum_{1\leq s_1 \leq t_1 < s_2 \leq t_2 <\cdots< s_j\leq t_j\leq n}\prod_{k=1}^j \omega_{2s_{k}-1}\omega_{2t_k},
\end{align}
where $\Pi_{k=1}^0\cdot$ denotes $1$ by convention.
  \item If $m=2n+1$ is odd, then for $0\leq j\leq n$ we have
\begin{align}\label{eq:rho-n}
\rho_{n-j}(\alpha)
&=\sum_{s=j}^n \rho_{s-j}(\omega_1,\ldots,\omega_{2s})\,\omega_{2s+1}\\
&=\sum_{1\le s_0 \le t_1 < s_1\le t_2 < s_2 \le \cdots \le t_j < s_j \le n+1}
\omega_{2s_0-1}\prod_{k=1}^j \omega_{2t_k}\,\omega_{2s_k-1}.
\end{align}
\end{enumerate}
\end{lemma}
\begin{lemma}
For $a,b \in \mathbb{F}_q^*$, let $\omega_{2s-1}^{\prime}=a\omega_{2s-1}$ and $\omega_{2s}^{\prime} =b\omega_{2s}$ for $s=1,2,\ldots$, then we have
\begin{gather}\label{eq:rho-a}
\rho_{n-j}(\omega_1^{\prime},\ldots,\omega_{2n}^{\prime})
=a^jb^j\rho_{n-j}(\omega_1,\ldots,\omega_{2n}),
\\
\label{eq:rho-b}
\rho_{n-j}(\omega_1^{\prime},\ldots,\omega_{2n+1}^{\prime})
=a^{j+1}b^j\rho_{n-j}(\omega_1,\ldots,\omega_{2n+1}).
\end{gather}
\end{lemma}

For nonnegative integers $m, n, k$ with $m\geq n\geq 2$, let $\Omega_{m,n,k}(\omega_1,\ldots,\omega_m)$ denote the following summation
\begin{gather}
\sum_i\big(\rho_{i}(\omega_1,\ldots,\omega_m)
\rho_{k-i}(\omega_2,\ldots,\omega_n)-\rho_{i}(\omega_2,\ldots,\omega_m)
\rho_{k-i}(\omega_1,\ldots,\omega_n)\big).
\end{gather}

\begin{lemma}\label{lem:Theta}
Let $m,n,k,$ be nonnegative integers with $m \geq n \geq 2$. Then,
\begin{align}\label{eq:Theta-rho}
&\Omega_{m+1,n+1,k+1}(\omega_1,\ldots,\omega_{m+1})=-\Omega_{m,n,k}(\omega_2,\ldots,\omega_{m+1}).
\end{align}
Furthermore, if $k< n \leq m$, then
\begin{align}\label{eq:Theta-m-n}
&\Omega_{m,n,k}(\omega_1,\ldots,\omega_{m})=0.
\end{align}
\end{lemma}
\noindent
{\bf Proof}: At first, $\eqref{eq:Theta-rho}$ follows from $\eqref{eq:rho1}$ directly.

If $k< n \leq m,$ according to $\eqref{eq:Theta-rho}$, without loss of generality, we assume $k=1< n\leq m$ further and then, from
\begin{align*}
\Omega_{m,n,k}(\omega_1,\ldots,\omega_m)=\rho_1(\omega_1,\ldots,\omega_m)\omega_2 \cdots \omega_n+\omega_1 \cdots \omega_m\rho_1(\omega_2,\ldots,\omega_n)\ \ \ \ \\
-\rho_1(\omega_2,\ldots,\omega_m)\omega_1 \cdots \omega_n-\omega_2 \cdots \omega_m\rho_1(\omega_1,\ldots,\omega_n)
\end{align*}
and $\eqref{eq:rho1}$, we see $\eqref{eq:Theta-m-n}.$
\eproof
\medskip

Since $\Lambda_{k,q}$ is edge-transitive, we may, without loss of generality, restrict attention to backtrackless walks starting from $[0]$ and $\langle 0\rangle$. Let $\Gamma=[l^{(1)}]\langle r^{(1)} \rangle [l^{(2)}] \langle r^{(2)} \rangle \cdots$ be such a given backtrackless walk of $\Lambda_{k,q}$, where $[l^{(1)}]=(0,0,\ldots,0)$ and $\langle r^{(1)} \rangle=(0,0,\ldots,0).$ For $i\geq 1$, let $x_i$ and $y_i$ denote the first entries (or colors) of $[l^{(i)}]$ and $\langle r^{(i)} \rangle$, respectively, and write
\begin{align}\label{eq:rho-x}
&\omega_{2i-1}=x_{i+1}-x_{i}, \quad \quad \omega_{2i}=y_{i+1}-y_{i}.
\end{align}
We note that $\omega_{j}$ is nonzero for all $j\geq1$. As a refinement of a closed-form expression given in \cite{Cheng14} for the backtrackless walks leading by $[l^{(1)}]=(0,0,\ldots,0)$, the following lemma was shown in \cite{XCT23}.

\begin{lemma}\label{lem:Gamma}
For any $i\geq 1$ and $j\geq 0$, let $l_{j}^{(i+1)}$ denote the $(j+1)$-th entries of $[l^{(i+1)}]$. Then, we have
\begin{align}\label{eq:rho-omega}
&l_{4j}^{(i+1)}=\rho_{i-j-1}(\omega_1,\ldots,\omega_{2i-1}),\\
&l_{4j+1}^{(i+1)}=\rho_{i-j-2}(\omega_2,\ldots,\omega_{2i-1}),\\
&l_{4j+2}^{(i+1)}=y_{i+1}l_{4j}^{(i+1)}-\rho_{i-j-1}(\omega_1,\ldots,\omega_{2i}),\\
&l_{4j+3}^{(i+1)}=y_{i+1}l_{4j+1}^{(i+1)}-\rho_{i-j-2}(\omega_2,\ldots,\omega_{2i}).
\end{align}
\end{lemma}

\medskip
For convenience, we say the walk $\Gamma$ is of type $(\omega_1,\omega_2,\ldots)$. If the first $2i$ vertices in the walk $\Gamma$ form a circuit in $\Lambda_{k,q}$ of length $2i$, we also call it a backtrackless circuit of type $(\omega_1,\ldots,\omega_{2i})$. We note that
\begin{align}\label{eq:rho-type-1-2i}
&\rho_{i-1}(\omega_2,\ldots,\omega_{2i})=\omega_2+\omega_4+\cdots+\omega_{2i}=y_{i+1}=y_{1}=0
\end{align}
is always a necessary condition for $\Lambda_{k,q}$ to have a backtrackless circuit of type $(\omega_1,\ldots,\omega_{2i})$.
Since the girth of $\Lambda_{k,q}$ is at least $k+4$, we have the following lemma.

\begin{lemma}\label{lem:backtrackless circuit}
If $\Lambda_{k,q}$ has a backtrackless circuit of length $2n$, then $k\leq 2n-4$.
\end{lemma}

\medskip

For any tuple $\alpha=(\omega_1,\omega_2,\ldots,\omega_{2n})$ of length $2n$, let
\begin{align*}
    \Delta_{2t-1}^{\alpha}=&\rho_{n-t}(\omega_2,\omega_3,\ldots,\omega_{2n}),\\
    \nabla_{2t-1}^{\alpha}=&\rho_{n-t}(\omega_1,\omega_2,\ldots,\omega_{2n-1}),\\
    \Delta_{2t}^{\alpha}=&\rho_{n-t-1}(\omega_2,\omega_3,\ldots,\omega_{2n-1}),\\
    \nabla_{2t}^{\alpha}=&\rho_{n-t}(\omega_1,\omega_2,\ldots,\omega_{2n}).
\end{align*}
Then, for $1\leq j\leq 2n-1$, according to Lemma \ref{lem:Theta} we see $\Omega_{2n,2n-1,2n-1-j}(\alpha)=0$ and thus from $$\Delta_{0}^{\alpha}=\rho_{n-1}(\omega_2,\omega_3,\ldots,\omega_{2n-1})=1$$
we have
\begin{align}\label{eq:nabla}
&\nabla_{2j}^{\alpha}=\sum_{t=0}^{j-1}(\nabla_{2t+1}^{\alpha}\Delta_{2j-2t-1}^{\alpha}-\nabla_{2t}^{\alpha}\Delta_{2j-2t}^{\alpha}).
\end{align}
According to Lemma~$\ref{lem:Gamma}$ and $\eqref{eq:nabla}$, one can deduce the following theorem easily (cf. \cite{XCT2023}).

\begin{theorem}\label{th:if and only if}
Let $\alpha=(\omega_1,\ldots,\omega_{2n})$ be a $2n$-tuple over $\mathbb{F}_q^*$.
\begin{enumerate}
  \item[1.] For $1\leq s \leq\lfloor (n-2)/2\rfloor,$ $\Lambda_{4s,q}$ has a backtrackless circuit of type $\alpha$ if and only if $\Delta_{i}^{\alpha}=\nabla_{i}^{\alpha}=0$ for $i=1,2,\ldots,2s+1$ .
  \item[2.] For $1\leq s \leq\lfloor (n-2)/2\rfloor,$ $\Lambda_{4s-1,q}$ has a backtrackless circuit of type $\alpha$ if and only if $\Delta_{2s+1}^{\alpha}=0$ and $\Delta_{i}^{\alpha}=\nabla_{i}^{\alpha}=0$ for $i=1,2,\ldots,2s$.
  \item[3.] For $k\in \{4s-2,4s-3\}$ with $1\leq s \leq\lfloor (n-1)/2\rfloor,$ $\Lambda_{k,q}$ has a backtrackless circuit of type $\alpha$ if and only if $\Delta_{i}^{\alpha}=\nabla_{i}^{\alpha}=0$ for $i=1,2,\ldots,2s$.
  \end{enumerate}
\end{theorem}

\section{Length-$(4s+4)$ Backtrackless Circuits in $\Lambda_{4s,q}$}\label{sec:backtrackless-Circuits}
Since for any graph $\Lambda_{4s,q}$ its girth is at least $4s+4$,
each of its length-$(4s+4)$ backtrackless circuits is a girth-cycle.
This section is devoted to proving the main theorem below, which further characterizes the necessary and sufficient conditions for the existence of length-$(4s+4)$ backtrackless circuits within $\Lambda_{4s,q}$.

\begin{theorem}\label{th:circuit}
Let $s\geq 1$. The graph $\Lambda_{4s,q}$ admits a backtrackless circuit of type $\alpha=(\omega_1,\ldots,\omega_{4s+4})$ if and only if
\begin{gather}\label{eq:rho-omega-ij}
\omega_k+\omega_{2s+2+k}=0,\quad 1\leq k \leq 2s+2,\\
\label{eq:rho-omega-ijs}
\rho_i(\omega_1,\ldots,\omega_{2s+2})=\rho_{i-1}(\omega_2,\ldots,\omega_{2s+1}), \quad 1\leq i \leq s.
\end{gather}
\end{theorem}

\medskip

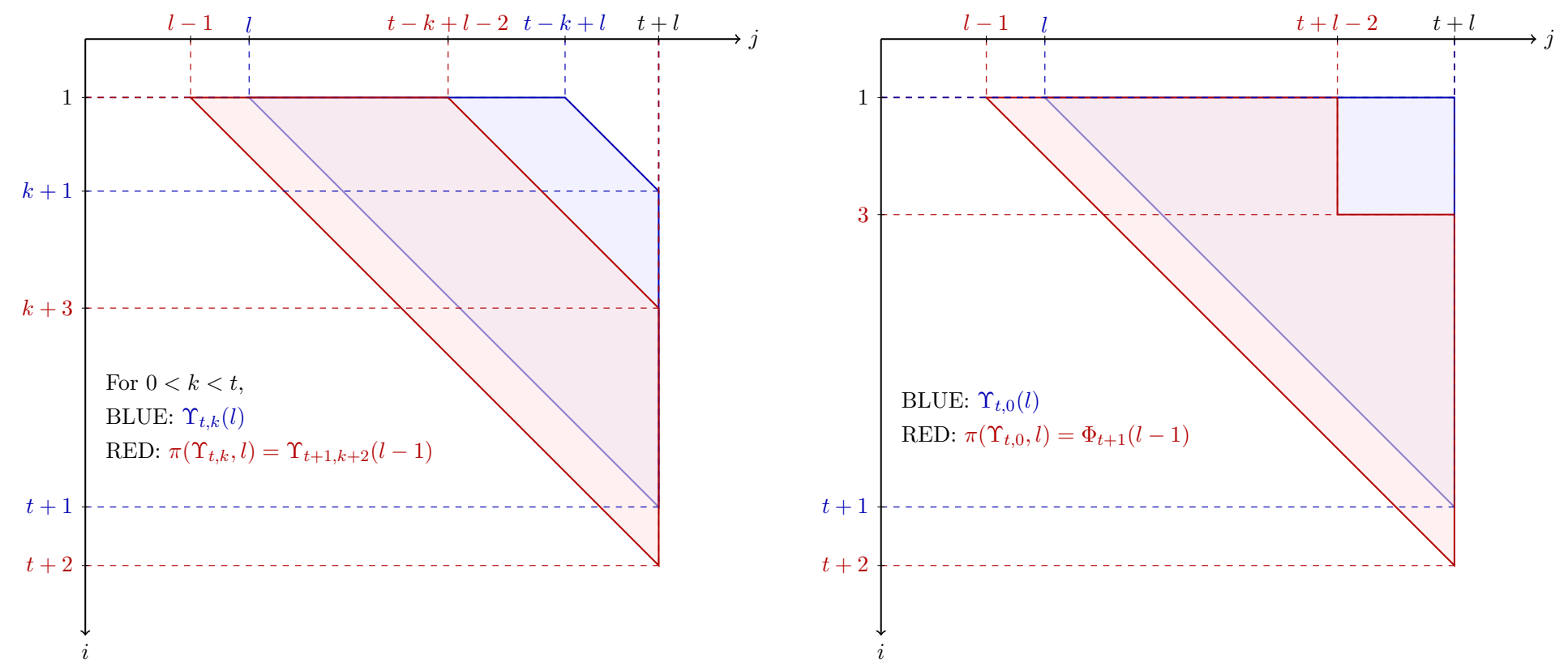
\begin{figure}[htbp]
\centering
\resizebox{\textwidth}{!}{%
\begin{tikzpicture}[
    x=1cm,
    y=-1cm,
    every node/.style={font=},
    axis/.style={->, thick},
    inregion/.style={
        draw=blue!70!black,
        fill=blue!10,
        fill opacity=0.55,
        line width=0.8pt
    },
    outregion/.style={
        draw=red!70!black,
        fill=red!10,
        fill opacity=0.55,
        line width=0.8pt
    },
    inproj/.style={
        blue!70!black,
        dashed,
        line width=0.55pt,
        opacity=0.85
    },
    outproj/.style={
        red!70!black,
        dashed,
        line width=0.55pt,
        opacity=0.85
    },
    ticklabel/.style={
        fill=white,
        fill opacity=1,
        text opacity=1,
        inner sep=1.2pt
    }
]


\newcommand{\axes}[1]{
    \draw[axis] (0,0) -- (11.2,0) node[right] {$j$};
    \draw[axis] (0,0) -- (0,10.2) node[below] {$i$};
    \node[align=left, anchor=west, ticklabel] at (0.3,6.5) {#1};
}

\newcommand{\xtickblue}[2]{
    \draw (#1,-0.06) -- (#1,0.06);
    \node[blue!70!black, above, yshift=2pt, ticklabel] at (#1,0) {$#2$};
}

\newcommand{\xtickred}[2]{
    \draw (#1,-0.06) -- (#1,0.06);
    \node[red!70!black, above, yshift=2pt, ticklabel] at (#1,0) {$#2$};
}

\newcommand{\xtickblack}[2]{
    \draw (#1,-0.06) -- (#1,0.06);
    \node[above, yshift=2pt, ticklabel] at (#1,0) {$#2$};
}

\newcommand{\ytickblue}[2]{
    \draw (-0.06,#1) -- (0.06,#1);
    \node[blue!70!black, anchor=east, ticklabel] at (-0.16,#1) {$#2$};
}

\newcommand{\ytickred}[2]{
    \draw (-0.06,#1) -- (0.06,#1);
    \node[red!70!black, anchor=east, ticklabel] at (-0.16,#1) {$#2$};
}

\newcommand{\ytickblack}[2]{
    \draw (-0.06,#1) -- (0.06,#1);
    \node[anchor=east, ticklabel] at (-0.16,#1) {$#2$};
}

\newcommand{\inprojection}[2]{
    \draw[inproj] (#1,0) -- (#1,#2);
    \draw[inproj] (0,#2) -- (#1,#2);
}

\newcommand{\outprojection}[2]{
    \draw[outproj] (#1,0) -- (#1,#2);
    \draw[outproj] (0,#2) -- (#1,#2);
}

\newcommand{\deletedholeprojections}[4]{

    \draw[outproj] (#1,0) -- (#1,#4);

    \draw[outproj] (0,#4) -- (#2,#4);

    \draw[inproj] (#2,0) -- (#2,#4);

    \draw[inproj] (0,#3) -- (#2,#3);
}


\begin{scope}

\coordinate (A) at (2.8,1);
\coordinate (B) at (8.2,1);
\coordinate (C) at (9.8,2.6);
\coordinate (D) at (9.8,8);

\coordinate (Aout) at (1.8,1);
\coordinate (Bout) at (6.2,1);
\coordinate (Cout) at (9.8,4.6);
\coordinate (Dout) at (9.8,9);

\path[inregion]  (A) -- (B) -- (C) -- (D) -- cycle;
\path[outregion] (Aout) -- (Bout) -- (Cout) -- (Dout) -- cycle;

\inprojection{2.8}{1}
\inprojection{8.2}{1}
\inprojection{9.8}{2.6}
\inprojection{9.8}{8}

\outprojection{1.8}{1}
\outprojection{6.2}{1}
\outprojection{9.8}{4.6}
\outprojection{9.8}{9}

\axes{
For $0<k<t$,\\[3pt]
BLUE: \textcolor{blue!70!black}{$\Upsilon_{t,k}(l)$}\\[3pt]
RED: \textcolor{red!70!black}{$\pi(\Upsilon_{t,k},l)=\Upsilon_{t+1,k+2}(l-1)$}
}

\xtickred{1.8}{l-1}
\xtickblue{2.8}{l}

\xtickred{6.2}{t-k+l-2}
\xtickblue{8.2}{t-k+l}

\xtickblack{9.8}{t+l}

\ytickblack{1}{1}
\ytickblue{2.6}{k+1}
\ytickred{4.6}{k+3}
\ytickblue{8}{t+1}
\ytickred{9}{t+2}

\end{scope}


\begin{scope}[xshift=13.6cm]

\coordinate (E) at (2.8,1);
\coordinate (F) at (9.8,1);
\coordinate (G) at (9.8,8);

\coordinate (H) at (1.8,1);
\coordinate (I) at (9.8,1);
\coordinate (J) at (9.8,9);

\coordinate (Q1) at (7.8,1);
\coordinate (Q2) at (9.8,1);
\coordinate (Q3) at (9.8,3);
\coordinate (Q4) at (7.8,3);

\path[inregion] (E) -- (F) -- (G) -- cycle;

\path[outregion] (H) -- (Q1) -- (Q4) -- (Q3) -- (J) -- cycle;

\inprojection{2.8}{1}
\inprojection{9.8}{1}
\inprojection{9.8}{8}

\outprojection{1.8}{1}
\outprojection{9.8}{9}

\deletedholeprojections{7.8}{9.8}{1}{3}

\axes{
BLUE: \textcolor{blue!70!black}{$\Upsilon_{t,0}(l)$}\\[3pt]
RED: \textcolor{red!70!black}{$\pi(\Upsilon_{t,0},l)=\Phi_{t+1}(l-1)$}
}

\xtickred{1.8}{l-1}
\xtickblue{2.8}{l}
\xtickred{7.8}{t+l-2}
\xtickblack{9.8}{t+l}

\ytickblack{1}{1}
\ytickred{3}{3}
\ytickblue{8}{t+1}
\ytickred{9}{t+2}

\end{scope}

\end{tikzpicture}
}
\caption{An illustration for Lemma \ref{lem:Psi-Phi}}
\label{fig-one}
\end{figure}

To prove this theorem, we introduce some notations and lemmas first.
For any subset $A\subset \mathbb{Z}\times \mathbb{Z}$ and integer $k\in\mathbb{Z}$, let $A(k)$ be the translate of $A$ defined by
\begin{align}\label{eq:A-k}
&A(k)=\{(i+1,j+k)|(i,j)\in A\},
\end{align}
and let $\pi(A,k)$ be the set of pairs $(i,j)$ for which either $\{(i,j+1),(i,j+2)\}\subseteq A(k)$ or $\{(i-1,j),(i-2,j)\}\subseteq A(k)$ is valid.
For $0\leq k\leq t$, let
\begin{gather}\label{a00}
\Upsilon_{t,k}=\{(i,j):0\leq i<k, 0\leq j-i\leq t-k\}
\cup\{(i,j): k\leq i\leq t, i\leq j\leq t\}.
\end{gather}
Clearly, $\Upsilon_{t,k}$ can also be expressed as
\begin{gather}\label{a01}
\{(i,j):0\leq j\leq t-k, 0\leq i\leq j\}
\cup\{(i,j): t-k<j\leq t, 0\leq j-i\leq t-k\}.
\end{gather}
Furthermore, from (\ref{a00}) and (\ref{a01}) we have the following lemma,
which is depicted in Figure~\ref{fig-one}.

\begin{lemma}\label{lem:Psi-Phi}
For $0\leq k<t$ and $l\in \mathbb{Z}$, we have
\begin{align}\label{eq:Pi-Upsilon-t}
&\pi(\Upsilon_{t,k},l)=\left\{
\begin{array}{ll}
\Upsilon_{t+1,\,k+2}(l-1),&\text{if }0<k<t,\\
\Phi_{t+1}(l-1),&\text{if }k=0,
\end{array}
\right.
\end{align}
where
\begin{align}
\label{b00}\Phi_{t+1}&=\Upsilon_{t+1,0}\setminus\{(0,t),(0,t+1),(1,t),(1,t+1)\}\\
&=\Upsilon_{t+1,1}\setminus\{(0,t),(1,t),(1,t+1)\}.
\label{b01}
\end{align}
\end{lemma}

\medskip

For any $B\subset \mathbb{Z}\times\mathbb{Z}$, let $G(B)$ denote the graph with $B$ as its vertex set such that two vertices $(i,j), (r,s)\in B$ are adjacent in $G(B)$ if and only if $\mid i-r \mid +\mid j-s \mid =1$. Clearly, for any $0\leq k< t$ the graph $G(\Upsilon_{t,k})$ is connected.
For the existence of length-$(4s+4)$ backtrackless circuits within $\Lambda_{4s,q}$, the following lemma claims a necessary and sufficient condition, which will be further refined in Theorem~\ref{th:circuit}.

\begin{lemma}\label{lem:type-alpha}
For $s\geq 1$, $\Lambda_{4s,q}$ has a backtrackless circuit of type $\alpha=(\omega_1,\ldots,\omega_{4s+4})$ if and only if
\begin{align}\label{eq:rho-k-s}
\rho_1(\omega_k,\ldots,\omega_{2s+2+k})=0\text{ for } 1\leq k \leq 2s+2,
\\
\label{eq:rho-k-s-j}
\rho_j(\omega_1,\ldots,\omega_{2s+2+j})=0\text{ for } 2\leq j \leq 2s+1.
\end{align}
\end{lemma}
\noindent{\bf Proof}:
We note that, according to Theorem~\ref{th:if and only if}, the graph $\Lambda_{4s,q}$ has a backtrackless circuit of type $\alpha=(\omega_1,\ldots,\omega_{4s+4})$ if and only if
\begin{gather}\label{a02}
\rho_k(\omega_2,\ldots,\omega_{4s+4})=0\text{ for } 2s+1\ge k \ge s+1,\\
\label{a03}\rho_k(\omega_1,\ldots,\omega_{4s+3})=0\text{ for } 2s+1\ge k \ge s+1,\\
\label{a04}\rho_k(\omega_2,\ldots,\omega_{4s+3})=0\text{ for } 2s\ge k \ge s+1,\\
\label{a05}\rho_k(\omega_1,\ldots,\omega_{4s+4})=0\text{ for } 2s+1\ge k \ge s+2.
\end{gather}

To prove the only-if-part of this lemma, we assume that (\ref{a02})-(\ref{a05}) are valid. Since $\rho_{2s+1}(\omega_i,\ldots,\omega_j)=0$ is valid for each
$(i,j)\in \{(1,4s+3),(1,4s+4),(2,4s+4)\}=\Upsilon_{1,0}(4s+3)$, from $\eqref{eq:rho1},\eqref{eq:rho2}$ and $\eqref{eq:Pi-Upsilon-t}$ we see $\rho_{2s}(\omega_i,\ldots,\omega_j)=0 $ for each $(i,j)\in \Phi_2(4s+2)=\{(1,4s+2),(3,4s+4)\}$. Then, according to (\ref{b00}) and (\ref{a02})-(\ref{a05}) we see $\rho_{2s}(\omega_i,\ldots,\omega_j)=0$ for each $(i,j)\in \Upsilon_{2,0}(4s+2)$. Furthermore, according to $\eqref{eq:rho1},\eqref{eq:rho2}$, $\eqref{eq:Pi-Upsilon-t}$, (\ref{b00}) and (\ref{a02})-(\ref{a05}), by induction one can show for $t=0,1,\ldots,s-1$ that
\begin{align}\label{eq:rho-2s+1-t}
&\rho_{2s+1-t}(\omega_i,\ldots,\omega_j)=0\text{ for each } (i,j)\in \Upsilon_{t+1,0}(4s+3-t).
\end{align}
Therefore, according to $\eqref{eq:rho1},\eqref{eq:rho2}$ and $\eqref{eq:Pi-Upsilon-t}$ we have $\rho_{s+1}(\omega_i,\ldots,\omega_j)=0$ for each $(i,j)\in \Phi_{s+1}(3s+3)$.
Then, from (\ref{b01}) and (\ref{a02})-(\ref{a04}) we see
$\rho_{s+1}(\omega_i,\ldots,\omega_j)=0$ for each $(i,j)\in \Upsilon_{s+1,1}(3s+3)$.
Hence, according to $\eqref{eq:rho1},\eqref{eq:rho2}$ and $\eqref{eq:Pi-Upsilon-t}$, by induction one can show for $t=0,1,\ldots,s$ that
\begin{align}\label{eq:rho-omega-i-j}
&\rho_{s+1-t}(\omega_i,\ldots,\omega_j)=0\text{ for each }(i,j)\in \Upsilon_{t+s+1,2t+1}(3s+3-t).
\end{align}
Since the set $\Upsilon_{t,k}$ contains the pair $(0,0)$, we see that $\eqref{eq:rho-k-s}$ and $\eqref{eq:rho-k-s-j}$ follow from $\eqref{eq:rho-2s+1-t}$, $\eqref{eq:rho-omega-i-j}$
and $\Upsilon_{2s+1,2s+1}(2s+3)=\{(i+1,i+2s+3):0\leq i\leq 2s+1\}$.

To prove the if-part of this lemma, we assume that $\eqref{eq:rho-k-s}$ and $\eqref{eq:rho-k-s-j}$ are valid. From $\eqref{eq:rho-k-s}$ we see $\eqref{eq:rho-omega-i-j}$ is valid for $t=s$. According to $\eqref{eq:rho2}$, $\rho_2(\omega_1,\ldots,\omega_{2s+4})=0$ and
$\rho_1(\omega_1,\ldots,\omega_{2s+3})=0$ we see $\rho_2(\omega_1,\ldots,\omega_{2s+5})=0$.
Then, according to $\eqref{eq:rho1}$ and $\rho_1(\omega_3,\ldots,\omega_{2s+5})=0$ we see
$\rho_2(\omega_2,\ldots,\omega_{2s+5})=0$.
Therefore, according to $\Upsilon_{2s+1,2s+1}(2s+3)=\pi(\Upsilon_{2s,2s-1},2s+4)$ and
the connectivity of $G(\Upsilon_{2s,2s-1})$ we see that $\eqref{eq:rho-omega-i-j}$ is valid for $t=s-1$.
Furthermore, according to $\eqref{eq:Pi-Upsilon-t}$ by induction one can show that $\eqref{eq:rho-omega-i-j}$ is valid for $t=s,s-1,\ldots,0$.
Since we have $\pi(\Upsilon_{s,0},3s+4)=\Phi_{s+1}(3s+3)\subset \Upsilon_{s+1,1}(3s+3)$, according to $\eqref{eq:rho1},\eqref{eq:rho2},\rho_{s+2}(\omega_1,\ldots,\omega_{3s+4})=0$ and the connectivity of $G(\Upsilon_{s,0})$ we see that $\eqref{eq:rho-2s+1-t}$ is valid for $t=s-1$. Furthermore, according to $\eqref{eq:rho-k-s-j}$ and
$$\pi(\Upsilon_{s-i,0},3s+4+i)=\Phi_{s+1-i}(3s+3+i)\subset \Upsilon_{s+1-i,\,0}(3s+3+i),1\leq i < s,$$
by induction one can show that $\eqref{eq:rho-2s+1-t}$ is valid for $t=s-1,s-2,\ldots,0$. Therefore, (\ref{a02})-(\ref{a05}) follow from $\eqref{eq:rho-2s+1-t}$,
\begin{gather*}
\{(2,4s+4),(1,4s+3),(2,4s+3)\}\subset \Upsilon_{s+1,1}(3s+3),\\
\{(2,4s+4),(1,4s+3),(1,4s+4)\}=\Upsilon_{1,0}(4s+3),
\end{gather*}
and
$$\{(2,4s+4),(1,4s+3),(2,4s+3),(1,4s+4)\}\subset \Upsilon_{t+1,0}(4s+3-t)$$
for $t=1,2,\ldots,s-1$.
\eproof
\medskip

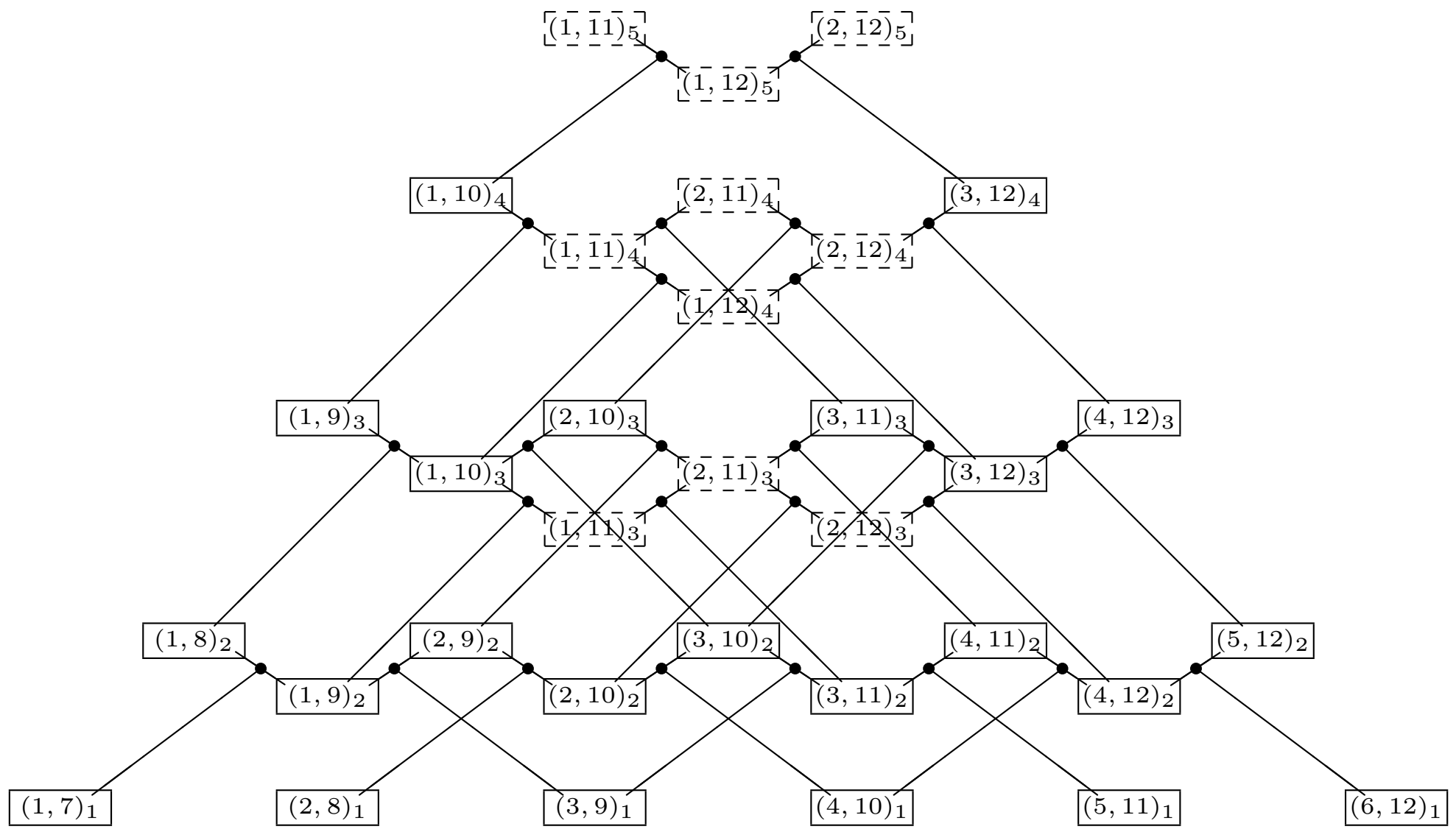
\begin{figure}[ht]
\newcounter{cms}
\setlength{\unitlength}{1mm}
\begin{picture}(160,75)
\label{ccc}
\put(48,70){\dashbox(9,3){\tiny $(1,11)_5$}}
\put(60,65){\dashbox(9,3){\tiny $(1,12)_5$}}
\put(72,70){\dashbox(9,3){\tiny $(2,12)_5$}}

\multiput(58.5,69)(12,0){2}{\circle*{1}}
\multiput(56.25,70.5)(24,0){1}{\line(3,-2){4,4}}
\multiput(68.25,67.5)(24,0){1}{\line(3,2){4,4}}
\multiput(58.5,69)(24,0){1}{\line(-4,-3){15,14}}
\multiput(70.5,69)(24,0){1}{\line(4,-3){15,14}}

\put(36,55){\framebox(9,3){\tiny $(1,10)_4$}}
\put(48,50){\dashbox(9,3){\tiny $(1,11)_4$}}
\put(60,45){\dashbox(9,3){\tiny $(1,12)_4$}}
\put(60,55){\dashbox(9,3){\tiny $(2,11)_4$}}
\put(72,50){\dashbox(9,3){\tiny $(2,12)_4$}}
\put(84,55){\framebox(9,3){\tiny $(3,12)_4$}}

\multiput(46.5,54)(12,0){4}{\circle*{1}}
\multiput(58.5,49)(12,0){2}{\circle*{1}}
\multiput(44.25,55.5)(24,0){2}{\line(3,-2){4,4}}
\multiput(56.25,52.5)(24,0){2}{\line(3,2){4,4}}
\multiput(56.25,50.5)(24,0){1}{\line(3,-2){4,4}}
\multiput(68.25,47.5)(24,0){1}{\line(3,2){4,4}}
\multiput(46.5,54)(24,0){2}{\line(-1,-1){16,16}}
\multiput(58.5,54)(24,0){2}{\line(1,-1){16,16}}
\multiput(58.5,49)(24,0){1}{\line(-1,-1){16,16}}
\multiput(70.5,49)(24,0){1}{\line(1,-1){16,16}}

\put(24,35){\framebox(9,3){\tiny $(1,9)_3$}}
\put(36,30){\framebox(9,3){\tiny $(1,10)_3$}}
\put(48,25){\dashbox(9,3){\tiny $(1,11)_3$}}
\put(48,35){\framebox(9,3){\tiny $(2,10)_3$}}
\put(60,30){\dashbox(9,3){\tiny $(2,11)_3$}}
\put(72,25){\dashbox(9,3){\tiny $(2,12)_3$}}
\put(72,35){\framebox(9,3){\tiny $(3,11)_3$}}
\put(84,30){\framebox(9,3){\tiny $(3,12)_3$}}
\put(96,35){\framebox(9,3){\tiny $(4,12)_3$}}

\multiput(34.5,34)(12,0){6}{\circle*{1}}
\multiput(46.5,29)(12,0){4}{\circle*{1}}
\multiput(32.25,35.5)(24,0){3}{\line(3,-2){4,4}}
\multiput(44.25,32.5)(24,0){3}{\line(3,2){4,4}}
\multiput(44.25,30.5)(24,0){2}{\line(3,-2){4,4}}
\multiput(56.25,27.5)(24,0){2}{\line(3,2){4,4}}
\multiput(34.5,34)(24,0){3}{\line(-1,-1){16,16}}
\multiput(46.5,34)(24,0){3}{\line(1,-1){16,16}}
\multiput(46.5,29)(24,0){2}{\line(-1,-1){16,16}}
\multiput(58.5,29)(24,0){2}{\line(1,-1){16,16}}

\put(12,15){\framebox(9,3){\tiny $(1,8)_2$}}
\put(24,10){\framebox(9,3){\tiny $(1,9)_2$}}
\put(36,15){\framebox(9,3){\tiny $(2,9)_2$}}
\put(48,10){\framebox(9,3){\tiny $(2,10)_2$}}
\put(60,15){\framebox(9,3){\tiny $(3,10)_2$}}
\put(72,10){\framebox(9,3){\tiny $(3,11)_2$}}
\put(84,15){\framebox(9,3){\tiny $(4,11)_2$}}
\put(96,10){\framebox(9,3){\tiny $(4,12)_2$}}
\put(108,15){\framebox(9,3){\tiny $(5,12)_2$}}

\multiput(22.5,14)(12,0){8}{\circle*{1}}
\multiput(20.25,15.5)(24,0){4}{\line(3,-2){4,4}}
\multiput(32.25,12.5)(24,0){4}{\line(3,2){4,4}}
\multiput(22.5,14)(24,0){4}{\line(-4,-3){15,14}}
\multiput(34.5,14)(24,0){4}{\line(4,-3){15,14}}

\put(0,0){\framebox(9,3){\tiny $(1,7)_1$}}
\put(24,0){\framebox(9,3){\tiny $(2,8)_1$}}
\put(48,0){\framebox(9,3){\tiny $(3,9)_1$}}
\put(72,0){\framebox(9,3){\tiny $(4,10)_1$}}
\put(96,0){\framebox(9,3){\tiny $(5,11)_1$}}
\put(120,0){\framebox(9,3){\tiny $(6,12)_1$}}

\end{picture}
\caption{\small Equations associated with dashed-outline boxes are algebraically equivalent to those for bottom-aligned boxes and leftmost boxes.}
\label{ccc}
\end{figure}

\medskip
Figure~\ref{ccc} provides a graphical interpretation of the proof of Lemma~\ref{lem:type-alpha} for the case $s=2$, where the notation $(i,j)_t$ denotes the equation $\rho_t(\omega_i,\ldots,\omega_j)=0$.
In this figure, the bottom-aligned boxes correspond equations $\eqref{eq:rho-k-s}$ while the leftmost boxes are associated with equations $\eqref{eq:rho-k-s-j}$, and dashed-outline boxes match equations (\ref{a02})-(\ref{a05}). Each encircled node interconnects three boxes, which allows the equation of any single box to be derived algebraically from the remaining two counterparts according to $\eqref{eq:rho1}$ and $\eqref{eq:rho2}$.
Suppose all equations for dashed-outline boxes hold. Then, for $t=1,2,\ldots, 2s$ the equations indexed $(\cdot,\cdot)_{2s+1-t}$ can be deduced from their counterparts $(\cdot,\cdot)_{2s+2-t}$.
Meanwhile, if equations for bottom-aligned boxes are valid, then for $t=1,2,\ldots, 2s$
the equations for $(\cdot,\cdot)_{t+1}$ can be derived from the equations indexed $(\cdot,\cdot)_{t}$ and the leftmost box $(1,2s+3+t)_{t+1}$.

\medskip

\noindent{\bf Proof of Theorem~\ref{th:circuit}}:
To prove the only-if-part, we assume that $\eqref{eq:rho-k-s}$ and $\eqref{eq:rho-k-s-j}$ are valid.
For $1\leq k \leq 2s+1$, according to $\eqref{eq:rho1}$ and $\eqref{eq:rho2}$ we have
\begin{align*}
&\omega_{2s+3+k}\rho_1(\omega_k,\ldots,\omega_{2s+2+k})-\omega_k\rho_1(\omega_{k+1},\ldots,\omega_{2s+3+k})\\
=&\rho_0(\omega_{k+2},\ldots,\omega_{2s+3+k})-\rho_0(\omega_k,\ldots,\omega_{2s+1+k})\\
=&(\omega_{2s+3+k}\omega_{2s+2+k}-\omega_k\omega_{k+1})\omega_{k+2}\cdots\omega_{2s+1+k}
\end{align*}
and then from $\eqref{eq:rho-k-s}$ we see
$$\omega_{2s+3+k}\omega_{2s+2+k}-\omega_k\omega_{k+1}=0,\quad 1 \leq k \leq 2s+1.$$
Hence, there is a nonzero element $\mu\in\mathbb{F}_q^*$ such that, for $1\leq i \leq s+1$,
\begin{align}\label{eq:frac}
&\frac{\omega_{2i-1}}{\omega_{2i+2s+1}}=\frac{\omega_{2i+2s+2}}{\omega_{2i}}=\mu.
\end{align}

Now we conclude, for $i\geq 1$, that
\begin{align}\label{eq:mu}
&\rho_i(\omega_1,\ldots,\omega_{2s+2})=-\mu\rho_{i-1}(\omega_2,\ldots,\omega_{2s+1}).
\end{align}
Indeed, from $\eqref{eq:rho2}$ and $\eqref{eq:frac}$ we have
\begin{align}
\rho_1(\omega_1,\ldots,\omega_{2s+3})=&\omega_{2s+3}\rho_1(\omega_1,\ldots,\omega_{2s+2})
+\rho_0(\omega_1,\ldots,\omega_{2s+1})\nonumber\\
=&\mu^{-1}\omega_1\big(\rho_1(\omega_1,\ldots,\omega_{2s+2})
+\mu\rho_0(\omega_2,\ldots,\omega_{2s+1})\big),\label{ab00}
\end{align}
and thus from $\eqref{eq:rho-k-s}$ we see
\begin{gather*}
\rho_1(\omega_1,\ldots,\omega_{2s+2})
=-\mu\rho_0(\omega_2,\ldots,\omega_{2s+1}),
\end{gather*}
i.e., $\eqref{eq:mu}$ is valid for $i=1$. From $\eqref{eq:rho1}$ and $\eqref{eq:frac}$ we have
\begin{align}
&\rho_2(\omega_1,\ldots,\omega_{2s+4})\nonumber\\
=&\omega_{2s+4}(\omega_{2s+3}\rho_2(\omega_1,\ldots,\omega_{2s+2})
+\rho_1(\omega_1,\ldots,\omega_{2s+1}))+\rho_1(\omega_1,\ldots,\omega_{2s+2})\nonumber\\
=&\omega_1\omega_2\rho_2(\omega_1,\ldots,\omega_{2s+2})
+\mu\big(\omega_2\rho_1(\omega_1,\ldots,\omega_{2s+1})-\rho_0(\omega_2,
\ldots,\omega_{2s+1})\big)\nonumber\\
=&\omega_1\omega_2\big(\rho_2(\omega_1,\ldots,\omega_{2s+2})
+\mu\rho_1(\omega_2,\ldots,\omega_{2s+1})\big),\label{ab01}
\end{align}
and thus from $\eqref{eq:rho-k-s-j}$ we see that $\eqref{eq:mu}$ is valid for $i=2$.
Now we assume $k\geq 1$ and $\eqref{eq:mu}$ is valid for $1\leq i \leq 2k$. From $\eqref{eq:rho-k-1-2},\eqref{eq:rho-a},\eqref{eq:rho-b}$ and $\eqref{eq:frac}$ we have
\begin{align}
&\rho_{2k+1}(\omega_1,\ldots,\omega_{2s+3+2k})\nonumber\\
=&\sum_{j}\Big(\rho_j(\omega_1,\ldots,\omega_{2s+2})
\rho_{2k+1-j}(\omega_{2s+3},\ldots,\omega_{2s+3+2k})\nonumber\\
&\ \ \ \ \ \ \ \ \ +\rho_j(\omega_1,\ldots,\omega_{2s+1})
\rho_{2k-j}(\omega_{2s+4},\ldots,\omega_{2s+3+2k})\Big)\nonumber\\
=&\sum_{j}\Big(\rho_j(\omega_1,\ldots,\omega_{2s+2})
\rho_{2k+1-j}(\omega_{1},\ldots,\omega_{2k+1})\mu^{-1}\nonumber\\
&\ \ \ \ \ \ \ \ \ +\rho_j(\omega_1,\ldots,\omega_{2s+1})
\rho_{2k-j}(\omega_{2},\ldots,\omega_{2k+1})\Big)\nonumber\\
=&\rho_{2k+1}(\omega_1,\ldots,\omega_{2s+2})
\rho_0(\omega_1,\ldots,\omega_{2k+1})\mu^{-1}\nonumber\\
&\ \ \ +\rho_{2k}(\omega_2,\ldots,\omega_{2s+1})
\rho_0(\omega_1,\ldots,\omega_{2k+1})\nonumber\\
&\ \ \ +\sum_{j}\Big(-\rho_{j-1}(\omega_2,\ldots,\omega_{2s+1})
\rho_{2k+1-j}(\omega_1,\ldots,\omega_{2k+1})\nonumber\\
&\ \ \ \ \ \ \ \ \ \ \ \ \ \ \ \ +\rho_j(\omega_1,\ldots,\omega_{2s+1})
\rho_{2k-j}(\omega_2,\ldots,\omega_{2k+1})\Big)\nonumber\\
=&\mu^{-1}\rho_0(\omega_1,\ldots,\omega_{2k+1})
\big(\rho_{2k+1}(\omega_1,\ldots,\omega_{2s+2})
+\mu\rho_{2k}(\omega_2,\ldots,\omega_{2s+1})\big)\nonumber\\
&\ \ \ +\Omega_{2s+1,2k+1,2k}(\omega_1,\ldots,\omega_{2s+1}),\label{ab02}
\end{align}
and thus from $\eqref{eq:Theta-m-n}$ and $\eqref{eq:rho-k-s-j}$ we see that $\eqref{eq:mu}$ is valid for $i=2k+1$. Furthermore, from $\eqref{eq:rho-k-1-2},\eqref{eq:rho-a},\eqref{eq:rho-b}$ and $\eqref{eq:frac}$ we have
\begin{align}
&\rho_{2k+2}(\omega_1,\ldots,\omega_{2s+4+2k})\nonumber\\
=&\sum_j\Big(\rho_j(\omega_1,\ldots,\omega_{2s+2})
\rho_{2k+2-j}(\omega_{2s+3},\ldots,\omega_{2s+4+2k})\nonumber\\
&\ \ \ \ \ \ +\rho_j(\omega_1,\ldots,\omega_{2s+1})
\rho_{2k+1-j}(\omega_{2s+4},\ldots,\omega_{2s+4+2k})\Big)\nonumber\\
=&\sum_j\Big(\rho_j(\omega_1,\ldots,\omega_{2s+2})
\rho_{2k+2-j}(\omega_{1},\ldots,\omega_{2k+2})\nonumber\\
&\ \ \ \ \ \ +\rho_j(\omega_1,\ldots,\omega_{2s+1})
\rho_{2k+1-j}(\omega_{2},\ldots,\omega_{2k+2})\mu\Big)\nonumber\\
=&\rho_{2k+2}(\omega_1,\ldots,\omega_{2s+2})\rho_0(\omega_1,\ldots,\omega_{2k+2})\nonumber\\
&\ \ \ +\mu\rho_{2k+1}(\omega_2,\ldots,\omega_{2s+1})\rho_0(\omega_1,\ldots,\omega_{2k+2})\nonumber\\
&\ \ \ +\mu\sum_{j}\Big(-\rho_{j-1}(\omega_2,\ldots,\omega_{2s+1})
\rho_{2k+2-j}(\omega_1,\ldots,\omega_{2k+2})\nonumber\\
&\ \ \ \ \ \ \ \ \ \ \ \ \ \ \ \ \ +\rho_j(\omega_1,\ldots,\omega_{2s+1})
\rho_{2k+1-j}(\omega_2,\ldots,\omega_{2k+2})\Big)\nonumber\\
=&\rho_0(\omega_1,\ldots,\omega_{2k+2})\big(\rho_{2k+2}(\omega_1,\ldots,\omega_{2s+2})
+\mu\rho_{2k+1}(\omega_2,\ldots,\omega_{2s+1})\big)\nonumber\\
&\ \ \ +\mu\Theta_{2s+1,2k+2,2k+1}(\omega_1,\ldots,\omega_{2s+1}),\label{ab03}
\end{align}
and thus from $\eqref{eq:Theta-m-n}$ and $\eqref{eq:rho-k-s-j}$ we see that $\eqref{eq:mu}$ is valid for $i=2k+2.$ Hence, by induction we see that \eqref{eq:mu} is valid for all $i\geq 1$.

Therefore, from $\rho_{s+1}(\omega_1,\ldots,\omega_{2s+2})=\rho_{s}(\omega_2,\ldots,\omega_{2s+1})=1 $ and $\eqref{eq:mu}$ we see $\mu=-1$, and then we have $\eqref{eq:rho-omega-ij}$ and $\eqref{eq:rho-omega-ijs}$.

To prove the if-part, we assume $\eqref{eq:rho-omega-ij}$ and $\eqref{eq:rho-omega-ijs}$ are valid.
We note that, by checking the proofs of (\ref{ab00})-(\ref{ab03}) accordingly, one can deduce $\eqref{eq:rho-k-s}$ and $\eqref{eq:rho-k-s-j}$ simply from $\eqref{eq:rho-omega-ij}$ and $\eqref{eq:rho-omega-ijs}$.

Therefore, according to Lemma~\ref{lem:type-alpha} we see the validity of this theorem.
\eproof

\medskip

\section{Explicit Girth-Cycles of $\Lambda_{4s,q}$}\label{sec:xtty}

Since each backtrackless circuit giving in Theorem~\ref{th:circuit} must be a girth-cycle of $\Lambda_{4s,q}$, which will be called a girth-cycle of pattern $(\omega_1,\omega_2,\ldots, \omega_{2s+2})$ hereafter.
By using of the necessary and sufficient condition given in Theorem~\ref{th:circuit}, in this section we will characterize all the girth-cycles of $\Lambda_{8,q}$, and all the graphs $\Lambda_{4s,q}$ which admit a girth-cycle of pattern $(x,1,\ldots,1,y)$.

\subsection{Characterization of Girth-Cycles in $\Lambda_{8,q}$}

By Theorem~\ref{th:circuit},
$\Lambda_{8,q}$ contains a girth-cycle of pattern $(\omega_1,\ldots,\omega_6)$ if and only if
\begin{gather}\label{8q46}
\left\{
\begin{array}{l}
(\omega_1-\omega_3-\omega_5)\omega_2+(\omega_1+\omega_3-\omega_5)\omega_4
+(\omega_1+\omega_3+\omega_5)\omega_6=0,\\
\omega_3(\omega_5-\omega_1)\omega_2\omega_4
-\omega_1(\omega_3+\omega_5)\omega_2\omega_6-\omega_5(\omega_1+\omega_3)\omega_4\omega_6=0.
\end{array}
\right.
\end{gather}
Substituting
\begin{gather}
a_1=\omega_1,\ \ a_2=-\omega_3,\ \ a_3=\omega_5,\ \ b_1=\omega_4,\ \ b_2=-\omega_6,\ \ b_3=-\omega_2,
\end{gather}
system (\ref{8q46}) can be rewritten in the following more symmetric form:
\begin{align}\label{8q50}
  \left\{
  \begin{array}{l}
  (a_1-a_2-a_3)b_1+(a_2-a_3-a_1)b_2+(a_3-a_1-a_2)b_3=0,\\
  a_1(a_2-a_3)b_2b_3+a_2(a_3-a_1)b_3b_1+a_3(a_1-a_2)b_1b_2=0.
  \end{array}
  \right.
\end{align}

Indeed, we have the following lemma which can be checked simply.
\begin{lemma}\label{aalem00}
For any elements $x_1,x_2,x_3,y_1,y_2,y_3\in\mathbb{F}_q^*$, the tuple $$(a_1,a_2,a_3,b_1,b_2,b_3)=(x_1,x_2,x_3,y_1,y_2,y_3)$$ is a solution of (\ref{8q50}) if and only if there exists a permutation $\sigma$ of $\{1,2,3\}$ such that $$(a_1,a_2,a_3,b_1,b_2,b_3)=
(y_{\sigma(1)},y_{\sigma(2)},y_{\sigma(3)},x_{\sigma(1)},x_{\sigma(2)},x_{\sigma(3)})$$
is a solution of (\ref{8q50}).
\end{lemma}

Then, all girth cycles of $\Lambda_{8,q}$ can be determined according to the following two lemmas, which can be derived from Lemma~\ref{aalem00} simply.
\begin{lemma}\label{aalem01}
If the characterize of $\mathbb{F}_q$ is $2$, then $\epsilon=(a_1,a_2,a_3,b_1,b_2,b_3)$ is a solution of (\ref{8q50}) if and only if $q\geq 4$ and there $c,d,e,f\in\mathbb{F}_q^*$ with $c\neq d$ and $e\neq f$ such that $\epsilon$ is one of the following two tuples
\begin{align}
(c,d,c+d,c^2e^{-1},d^2f^{-1},(c+d)^2(e+f)^{-1}),\label{pp00}\\
(c^2e^{-1},d^2f^{-1},(c+d)^2(e+f)^{-1},c,d,c+d).
\end{align}
\end{lemma}
\noindent{\bf Proof}: The if-part can be checked easily. To prove the only-if-part, we assume $\epsilon=(a_1,a_2,a_3,b_1,b_2,b_3)$ is a solution of (\ref{8q50}).
Since the first equation in (\ref{8q50}) is equivalent to $(a_1+a_2+a_3)(b_1+b_2+b_3)=0$,
according to Lemma~\ref{aalem00} we can assume $a_1+a_2+a_3=0$ without loss of generality. Let $c=a_1$ and $d=a_2$. Then, we have
$c\neq d$, $a_3=c+d$ and thus from the second equation in (\ref{8q50}) we see
$$
c^2b_1^{-1}+d^2b_2^{-1}+(c+d)^2b_3^{-1}=0.
$$
For $e=c^2b_1^{-1}$ and $f=d^2b_2^{-1}$, we have $e\neq f$, $(c+d)^2b_3^{-1}=e+f$ and (\ref{pp00}).
\eproof

\begin{lemma}\label{aalem02}
Let $\epsilon=(a_1,a_2,a_3,b_1,b_2,b_3)$ be an arbitrary tuple over $\mathbb{F}_q^*$, where $q$ is a power of odd prime $p$. Then, we have the following conclusions:
\begin{enumerate}
\item If $a_1=a_2=a_3$, then $\epsilon$ is a solution of (\ref{8q50}) if and only if $b_1+b_2+b_3=0$.
\item If $|\{a_1,a_2,a_3\}|=2$, then $\epsilon$ is a solution of (\ref{8q50}) if and only if $q\geq 5$ and there exist $c,d,e\in\mathbb{F}_q^*$ with $c\not\in\{d,2d\}$ such that $\epsilon$ is one of the following three tuples
  \begin{gather}
  (c,d,d,2ce,(c-2d)e,(c-2d)e),\label{pp30}\\
  (d,c,d,(c-2d)e,2ce,(c-2d)e),\label{pp31}\\
  (d,d,c,(c-2d)e,(c-2d)e,2ce).\label{pp32}
  \end{gather}
\item Suppose $a_1,a_2,a_3$ are distinct. Let
\begin{align*}
&\Theta(a_1,a_2,a_3)=(a_1+a_2-a_3)(a_2+a_3-a_1)(a_3+a_1-a_2).
\end{align*}
\begin{enumerate}
\item If $\Theta(a_1,a_2,a_3)=0$, then $\epsilon$ is a solution of (\ref{8q50}) if and only if $q\geq 5$ and there exist $c,d,e\in\mathbb{F}_q^*$ with $c\not\in\{d,-d\}$ such that $\epsilon$ is one of the following three tuples
\begin{gather}
(c,d,c+d,ce,-de,(c-d)e),\label{pp40}\\
(d,c+d,c,-de,(c-d)e,ce),\\
(c+d,c,d,(c-d)e,ce,-de).
\end{gather}
\item If $\Theta(a_1,a_2,a_3)\neq 0$, then $\epsilon$ is a solution of (\ref{8q50}) if and only if $q\geq 7$ and there exists $d\in\mathbb{F}_q$ such that
  \begin{align}
d^2 &=a_1a_2a_3(a_1a_2a_3-\Theta(a_1,a_2,a_3)),\label{8q59}\\  b_1&=\frac{(d-a_1a_2(a_1-a_2))b_2}{a_2(a_3-a_1)(a_1-a_2-a_3)},\label{8q60}\\
  b_3&=\frac{(d-a_3a_2(a_3-a_2))b_2}{a_2(a_1-a_3)(a_3-a_2-a_1)}.\label{8q61}
  \end{align}
\end{enumerate}
\end{enumerate}
\end{lemma}
\noindent{\bf Proof}: Since the if-part of each conclusion can be checked easily, we only give proofs for the only-if-parts.

The only-if-part of the conclusion 1 is obvious.

To show the only-if-part of the conclusion 2, according to Lemma~\ref{aalem00} without loss of generality we assume $a_1=c\neq a_2=a_3=d$.
From the second equation of (\ref{8q50}) we see $b_2=b_3$.
Therefore, from the first equation of (\ref{8q50}) we see $(c-2d)b_1=2cb_2$, and then
we have $c\neq 2d$, $b_1=2ce$ and $b_2=b_3=(c-2d)e$ for some $e\in\mathbb{F}_q^*$. Hence, we have $q\geq 5$ and (\ref{pp30}).

To show the only-if-part of the conclusion 3(a), according to Lemma~\ref{aalem00} without loss of generality we assume $a_1=c,a_2=d,a_3=c+d$ are distinct. Clearly, we have $c\not\in\{d,-d\}$ and $q\geq 5$. Since from the first equation of (\ref{8q50}) we have $db_1+cb_2=0$,
we see $b_1=ce$ and $b_2=-de$ for some $e\in\mathbb{F}_q^*$. Therefore, from the second equation of (\ref{8q50}) we see $b_3=(c-d)e$. Hence, we have (\ref{pp40}).

To show the only-if-part of the conclusion 3(b), we assume $a_1,a_2,a_3$ are distinct and $\Theta(a_1,a_2,a_3)\neq 0$. Clearly, there are $b_1,b_2,b_3\in\mathbb{F}_q^*$ satisfying (\ref{8q50})
if and only if there are $b_1,b_2\in\mathbb{F}_q^*$ such that
\begin{align*}
&(a_3-a_1-a_2)a_3(a_1-a_2)b_1b_2\\
=&\big((a_1-a_2-a_3)b_1+(a_2-a_3-a_1)b_2\big)\big(a_1(a_2-a_3)b_2+a_2(a_3-a_1)b_1\big),
\end{align*}
which is equivalent to
\begin{align*}
a_2(a_3-a_1)(a_1-a_2-a_3)b_1^2+a_1(a_2-a_3)(a_2-a_3-a_1)b_2^2+2a_1a_2(a_1-a_2)b_1b_2=0,
\end{align*}
i.e.,
\begin{align*}
&\big(a_2(a_3-a_1)(a_1-a_2-a_3)b_1b_2^{-1}+a_1a_2(a_1-a_2)\big)^2\nonumber\\
=&a_1^2a_2^2(a_1-a_2)^2-a_1a_2(a_3-a_1)(a_2-a_3)(a_1-a_2-a_3)(a_2-a_3-a_1)\nonumber\\
=&a_1^2a_2^2(a_1-a_2)^2-a_1a_2(a_3-a_1)(a_3-a_2)\big((a_1-a_2)^2-a_3^2\big)\nonumber\\
=&a_1a_2a_3\big(a_1a_2a_3+(a_1+a_2-a_3)((a_1-a_2)^2-a_3^2)\big)\nonumber\\
=&a_1a_2a_3(a_1a_2a_3-\Theta(a_1,a_2,a_3))
\end{align*}
is a quare in $\mathbb{F}_q$, and then for some $d$ with (\ref{8q59}) we have (\ref{8q60}) and
\begin{align*}
b_3=&\frac{(a_1-a_2-a_3)b_1+(a_2-a_3-a_1)b_2}{a_1+a_2-a_3}\\
=&\frac{(d-a_1a_2(a_1-a_2))b_2+a_2(a_3-a_1)(a_2-a_3-a_1)b_2}{a_2(a_3-a_1)(a_1+a_2-a_3)}\\
=&\frac{(d-a_3a_2(a_3-a_2))b_2}{a_2(a_1-a_3)(a_3-a_2-a_1)}.
\end{align*}
Hence, (\ref{8q61}) is also valid.
\eproof

\medskip
Here are some examples of Lemmas~\ref{aalem01} and \ref{aalem02}.
\begin{enumerate}
\item If $q=3$, then the equation (\ref{8q50}) has 4 solutions $(a_1,a_2,a_3,b_1,b_2,b_3)$:
\begin{gather*}
(1,1,1,1,1,1),\
(1,1,1,2,2,2),\
(2,2,2,1,1,1),\
(2,2,2,2,2,2),
\end{gather*}
which coincide those corresponding the girth-cycles of $\Lambda_{8,3}$ given in \cite{XCT23}.

\item If $q=4$, then the equation (\ref{8q50}) has 54 solutions $(a_1,a_2,a_3,b_1,b_2,b_3)$:
\begin{gather*}
(a,b,c,da,db,dc),\ (a,b,c,d,d,d),\ (d,d,d,a,b,c),
\end{gather*}
where $d\in\{a,b,c\}=\mathbb{F}_4^*$.

\item If $q=5$, then the equation (\ref{8q50}) has 144 solutions $(a_1,a_2,a_3,b_1,b_2,b_3)$:
\begin{gather*}
(a,a,a,b,b,3b),(a,a,a,b,3b,b),(a,a,a,3b,b,b),\\
(a,a,3a,b,b,b),(a,3a,a,b,b,b),(3a,a,a,b,b,b),\\
(a,a,4a,b,b,4b),(a,4a,a,b,4b,b),(4a,a,a,4b,b,b),
\end{gather*}
where $a,b\in\mathbb{F}_5^*=\{1,2,3,4\}$.

\item
If $q=7$, then the equation (\ref{8q50}) has 1188 solutions $(a_1,a_2,a_3,b_1,b_2,b_3)$:
\begin{gather*}
(a,a,a,b,2b,4b),\ (a,a,a,b,4b,2b),\\
(a,a,a,b,b,5b),\ (a,a,a,b,5b,b),\ (a,a,a,5b,b,b),\\
(a,a,3a,b,b,6b),\ (a,3a,a,b,6b,b),\ (3a,a,a,6b,b,b),\\
(a,a,4a,b,b,4b),\ (a,4a,a,b,4b,b),\ (4a,a,a,4b,b,b),\\
(a,a,5a,b,b,b),\ (a,5a,a,b,b,b),\ (5a,a,a,b,b,b),\\
(a,a,6a,b,b,3b),\ (a,6a,a,b,3b,b),\ (6a,a,a,3b,b,b),\\
(a,2a,6a,b,3b,5b),\ (a,3a,5a,b,2b,6b),\ (a,6a,2a,b,5b,3b),\ (a,5a,3a,b,6b,2b),\\
(2a,a,6a,3b,b,5b),\ (3a,a,5a,2b,b,6b),\ (6a,a,2a,5b,b,3b),\ (5a,a,3a,6b,b,2b),\\
(2a,6a,a,3b,5b,b),\ (3a,5a,a,2b,6b,b),\ (6a,2a,a,5b,3b,b),\ (5a,3a,a,6b,2b,b),\\
(a,2a,4a,b,b,b),\ (a,4a,2a,b,b,b),\\
(a,2a,4a,b,2b,4b),(a,4a,2a,b,4b,2b),
\end{gather*}
where $a,b\in\mathbb{F}_7^*=\{1,2,3,4,5,6\}$.
\end{enumerate}

\subsection{Girth-Cycles of Pattern $(x,1,\ldots,1,y)$ in $D(4s,q)$}

For $0\leq b\leq a$, let $C_a^b$ denote the binomial coefficient $\frac{a!}{\,b!(a-b)!\,}$.
\begin{lemma}
\label{aalem50}
Let $\varpi=(\omega_1,\omega_2,\ldots,\omega_{2s+2})$ be a $(2s+2)$-tuple over $\mathbb{F}_q^*$, where $\omega_1=x$, $\omega_{2s+2}=y$ and $\omega_i=1$ for $i=2,3,\ldots,2s+1$.
Then, the graph $\Lambda_{4s,q}$ has a backtrackless circuit of pattern $\varpi$ if and only if
\begin{gather}
(2s-1)xy=1-x-y ,\label{bb09}\\
C_{2s-k}^k\equiv(2s-1)C_{2s-k}^{k-1}(\bmod\ p), \text{ for }2\leq k\leq s,
\label{bb10}
\end{gather}
where $p$ is the characteristic of $\mathbb{F}_q$.
\end{lemma}
\noindent{\bf Proof}:
Clearly, we have $\rho_k(\omega_i,\ldots,\omega_j)=C_{j-i+1-k}^k$
for $2\leq i\leq j\leq 2s+1$.
Therefore, according to Theorem~\ref{th:circuit} and
\begin{align*}
\rho_{k-1}(\omega_2,\ldots,\omega_{2s+1})=&C_{2s-k+1}^{k-1}=C_{2s-k}^{k-1}+C_{2s-k}^{k-2},\\
\rho_k(\omega_1,\ldots,\omega_{2s+2})
=&xy\rho_k(\omega_2,\ldots,\omega_{2s+1})+x\rho_{k-1}(\omega_2,\ldots,\omega_{2s})\nonumber\\
&\ \ \ \ \ \ \ \ +y\rho_{k-1}(\omega_3,\ldots,\omega_{2s+1})+\rho_{k-2}(\omega_3,\ldots,\omega_{2s})\nonumber\\
=&xy C_{2s-k}^k+(x+y)C_{2s-k}^{k-1}+C_{2s-k}^{k-2},
\end{align*}
the graph $\Lambda_{4s,q}$ has a backtrackless circuit of pattern $\varpi$ if and only if
\begin{align*}
xy C_{2s-k}^k=(1-x-y)C_{2s-k}^{k-1},\text{ for }1\leq k\leq s,
\end{align*}
which is equivalent to (\ref{bb09}) and (\ref{bb10}).
\eproof

\medskip
According to Lucas' theorem, for integers
\begin{gather*}
a=a_0+a_1p+\cdots+a_tp^t, \\
b=b_0+b_1p+\cdots+b_tp^t,
\end{gather*}
where $a_i,b_i\in\{0,1,\ldots,p-1\}$, $i=0,1,\ldots,t$, we have
\begin{gather}
C_a^b\equiv C_{a_0}^{b_0} C_{a_1}^{b_1} \cdots C_{a_t}^{b_t}(\bmod\ p).
\label{bb07}
\end{gather}
By using (\ref{bb07}), we can deduce
the following lemma, which exhibits a simple criterion for (\ref{bb09}) and (\ref{bb10}).

\begin{lemma}\label{aalem51}
For any positive integer $s$, there are $x,y\in\mathbb{F}_q^*$ such that (\ref{bb09}) and (\ref{bb10}) if and only if there exist a positive integer $r$ such that
one of the following conditions holds
\begin{enumerate}
  \item $p=2$, $s=2^r-1$ and $1\in\{x,y\}$,
  \item $p\neq 2$, $q\neq 3$, $s=(p^r-1)/2$, $y\not\in \{1,2^{-1}\}$ and $x=2^{-1}(1-(2y-1)^{-1})$.
\end{enumerate}
\end{lemma}
\noindent{\bf Proof}: To show the only-if-part, we assume (\ref{bb09}) and (\ref{bb10}) hold.
From (\ref{bb10}) we get $$1=C_s^s\equiv(2s-1)C_s^{s-1}=(2s-1)s(\bmod\ p),$$
which is equivalent to
\begin{gather}
p|(2s+1)(s-1).
\label{bb11}
\end{gather}

At first we assume $p\nmid(2s+1)$. Then, we have $p|(s-1)$ and $s\geq p+1$. From $p\nmid 3s=(2s+1)+(s-1)$ we see $p\neq 3$. From (\ref{bb10}) and $p|(s-1)$ we get $$0\equiv 3\big((2s-1)C_{2s-3}^2 -C_{2s-3}^3\big)=
(2s-3)(s-2)(4s+2)\equiv 6(\bmod\ p),$$
and then we have $p=2$. Let $s+1=2^rd$, where $r\geq 1$ and $2\nmid d$. We now conclude $d=1$.
Assume $d>1$ in contrast. Then $2^{r+1}<s$ and $d=2e+3$ for some $e\geq 0$. According to (\ref{bb07}) and
$$2s-2^{r+1}=2(2^rd-1)-2^{r+1}=2^{r+2}e+\sum_{i=1}^{r+1}2^i$$
we see $C_{2s-2^{r+1}}^{2^{r+1}}\equiv 1(\bmod\ 2)$,
$C_{2s-2^{r+1}}^{2^{r+1}-1}\equiv 0(\bmod\ 2)$ and then for $k=2^{r+1}$ we have
\begin{gather*}
C_{2s-k}^k\not\equiv(2s-1)C_{2s-k}^{k-1}(\bmod\ 2),
\end{gather*}
which contradicts to (\ref{bb10}). Hence $d=1$ and $s=2^r-1$. According to $p=2$ and (\ref{bb09}) we see $1\in\{x,y\}$ easily. Therefore, the first condition is valid.

Now we assume $p|(2s+1)$. Let $2s+1=p^rd$, where $r\geq 1$ and $p\nmid d$.
We conclude $d=1$.
Assume in contrast that $d>1$. Then, from $2\nmid d$ we see $d\geq 3$ and thus $p^r<s$. Let $d-2=fp+e$, where $0\leq e<p$.
According to (\ref{bb07}) and
$$2s-p^r=p^r(d-2)+(p^r-1)=fp^{r+1}+ep^r+\sum_{i=0}^{r-1}(p-1)p^i$$
we see
$C_{2s-p^r}^{p^r}\equiv e(\bmod\ p)$,
$C_{2s-p^r}^{p^r-1}\equiv 1(\bmod\ p)$
and thus from (\ref{bb10}) we have
\begin{align*}
0\equiv C_{2s-p^r}^{p^r}-(2s-1)C_{2s-p^r}^{p^r-1}
\equiv& e+2\equiv d(\bmod\ p),
\end{align*}
which contradicts $p\nmid d$. Hence, we must have $d=1$ and
$s=(p^r-1)/2$.
Furthermore, from (\ref{bb09}) we see
\begin{align}
2xy-x-y+1=0
\end{align}
i.e., $(2x-1)(2y-1)=-1$. Therefore, $y\not\in\{1,2^{-1}\}$ and $x=2^{-1}(1-(2y-1)^{-1})$. Hence, $q>3$ and
then the second condition is valid.

To show the if-part, we assume one of the two conditions is valid.

At first, we assume the first condition is valid. Clearly, we have (\ref{bb09}).
From $2s=2^{r+1}-2=\sum_{i=1}^r2^i$ and (\ref{bb07}), for $1\leq k\leq s=2^r-1=\sum_{i=0}^{r-1}2^i$ we have
\begin{align*}
&C_{2s-k}^k
\equiv C_{2s-k}^{k-1}\equiv
\Big\{
\begin{array}{ll}
1(\bmod\ 2), &\text{if }k=2^j-1, 1\leq j\leq r,\\
0(\bmod\ 2), &\text{otherwise}.
\end{array}
\end{align*}
Hence, we have (\ref{bb10}).

Now we assume the second condition is valid. Clearly, we have (\ref{bb09}). Let $a$ be an arbitrary integer with $0\leq a\leq p^r$.
Since in $\mathbb{F}_p[[z]]$, the ring of formal power series, we have
\begin{align*}
\sum_{0\leq b\leq p^r-a}C_{p^r-a}^bz^b
=(1+z)^{p^r}(1+z)^{-a}
=(1+z^{p^r})\sum_{b\geq 0}(-1)^bC_{a+b-1}^bz^b,
\end{align*}
for any integer $b$ with $0\leq b\leq p^r-a$ we see
\begin{align}
C_{p^r-a}^b\equiv(-1)^bC_{a+b-1}^b(\bmod\ p).
\label{bb12}
\end{align}
Therefore, from $C_{2k}^k=2C_{2k-1}^{k-1}$ and (\ref{bb12}) we see (\ref{bb10}).
\eproof

\medskip
The following corollary follows directly from Lemmas~\ref{aalem50} and \ref{aalem51}.
\begin{corollary}\label{cor00}
Let $r, m$ be positive integers. Then,
\begin{enumerate}
  \item the girth of $\Lambda_{4(2^r-1),2^m}$ equals $2^{r+2}$;
  \item for odd prime $p$ with $p^m>3$, the girth of $\Lambda_{2(p^r-1),p^m}$ equals $2(p^r+1)$.
\end{enumerate}
\end{corollary}

Note that the following lemma is from \cite{Sun17} (cf. \cite{CZTX26}).
\begin{lemma}\label{lem70}
If $g(D(n,q))$ has girth $g_{n}$, then, for any $t\geq 0$ and
\begin{align*}
k=\left\{
\begin{array}{ll}
  g_np^t-5, & \text{if }n\equiv 3(\bmod 4), \\
  g_np^t-4, & \text{otherwise},
\end{array}
\right.
\end{align*}
the graph $D(k,q)$ has girth $g_k=g_np^t$.
\end{lemma}

Therefore, according to $D(k,q)\cong\Lambda_{k,q}$, Corollary~\ref{cor00} and Lemma~\ref{lem70},
we have the following theorem which confirms Conjecture~\ref{conmain} for the new parameter family
$k=2(p^r+1)p^s-4$
under the conditions $p>2$ and $q>3$.
\begin{theorem}\label{theo50}
Assume $p$ is an odd prime and $p^m>3$. Then, for any nonnegative integers $r,t$, the girth of $D(2(p^r+1)p^t-4,p^m)$ equals $2(p^r+1)p^t$.
\end{theorem}

\section{Concluding Remarks}

This paper investigates the girth cycles of the even-indexed subclass of Lazebnik-Ustimenko bipartite graphs $D(4s,q)$. Adopting the standard analytical framework established in \cite{Cheng14,Cheng16,XCT2023}, we conduct our analysis based on the isomorphic graph model $\Lambda_{k,q}$, which effectively simplifies the intricate algebraic structure of the original $D(k,q)$ construction. We establish a precise explicit necessary and sufficient condition for the existence of length-$(4s+4)$ backtrackless circuits in $\Lambda_{4s,q}$. Equipped with this newly derived criterion, we construct concrete girth cycles, achieve a complete characterization of all girth cycles of $\Lambda_{8,q}$, and classify all $\Lambda_{4s,q}$ graphs that possess girth cycles with the symmetric pattern $(x,1,\dots,1,y)$. Most importantly, our results validate Conjecture~\ref{conmain} for the previously unexplored parameter family $k=2(p^r+1)p^t-4$ over finite fields of odd characteristic $p>2$ with $q=p^m>3$, further extending the valid parameter range of this classic girth conjecture.

\renewcommand\refname{Reference}


\begin{thebibliography}{10}

\bibitem{Wenger91}
R. Wenger, Extremal graphs with no $C^{4s}, C^{6s}$, or $C^{10s}$,
J. Combin. Theory Ser. B 52(1) (1991) 113--116.

\bibitem{Lazebnik93}
F. Lazebnik, V. A. Ustimenko, New examples of graphs without small cycles and of large size,
European J. Combin. 14 (1993) 445--460.

\bibitem{LazebnikUstimenko95}
F. Lazebnik, V. A. Ustimenko, Explicit construction of graphs with an arbitrary large girth and of large size,
Discrete Appl. Math. 60 (1995) 275--284.

\bibitem{Lazebnik95}
F. Lazebnik, V. A. Ustimenko, A. J. Woldar, A new series of dense graphs of high girth,
Bull. Amer. Math. Soc. 32(1) (1995) 73--79.

\bibitem{Furedi95}
Z. F\"{u}redi, F. Lazebnik, A. Seress, V. A. Ustimenko, A. J. Woldar,
Graphs of prescribed girth and bi-degree,
J. Combin. Theory Ser. B 64(2) (1995) 228--239.

\bibitem{Lazebnik96}
F. Lazebnik, V. A. Ustimenko, A. J. Woldar,
A characterization of the component of the graphs $D(k,q)$,
Discrete Math. 157 (1996) 271--283.

\bibitem{Lazebnik97}
F. Lazebnik, V. A. Ustimenko, A. J. Woldar,
New upper bounds on the order of cages,
Electron. J. Combin. 14 (R13) (1997) 1--11.

\bibitem{Lazebnik01}
F. Lazebnik, A. J. Woldar,
General properties of families of graphs defined by some systems of equations,
J. Graph Theory 38(2) (2001) 65--86.

\bibitem{Lazebnik02}
F. Lazebnik, R. Viglione,
An infinite series of regular edge- but not vertex-transitive graphs,
J. Graph Theory 41(4) (2002) 249--258.

\bibitem{Ustimenko03}
V. A. Ustimenko, A. J. Woldar,
Extremal properties of regular and affine generalized $m$-gons as tactical configurations,
European J. Combin. 24 (2003) 99--111.

\bibitem{Lazebnik04}
F. Lazebnik, R. Viglione,
On the connectivity of certain graphs of high girth,
Discrete Math. 277 (2004) 309--319.

\bibitem{Kim04}
J. Kim, U. Peled, I. Perepelitsa, V. Pless, S. Friedland,
Explicit construction of families of LDPC codes with no 4-cycles,
IEEE Trans. Inform. Theory 50(10) (2004) 2378--2388.

\bibitem{Dmytrenko07}
V. Dmytrenko, F. Lazebnik, J. Williford,
On monomial graphs of girth eight,
Finite Fields Appl. 13 (2007) 828--842.

\bibitem{Futorny07}
V. Futorny, V. Ustimenko,
On small world semiplanes with generalised Schubert cells,
Acta Applicandae Math. 98 (2007) 47--61.

\bibitem{Ustimenko}
V. A. Ustimenko,
On linguistic dynamical systems, families of graphs of large girth, and cryptography,
J. Math. Sci. 140(3) (2007) 461--471.

\bibitem{Ustimenko09}
V. A. Ustimenko,
On the homogeneous algebraic graphs of large girth and their applications,
Linear Algebra Appl. 430(7) (2009) 1826--1837.

\bibitem{Yan11}
T. Yan, Y. Tang,
Constructions of LDPC codes based on polarity graphs with prescribed girth,
in: Asia-Pacific Youth Conference on Communication, 2011 APYCC, 2011, pp. 60--62.

\bibitem{Karimi}
M. Karimi, A. H. Banihashemi.
Efficient algorithm for finding dominant trap-\allowbreak ping sets of LDPC codes[J].
IEEE Trans. Inf. Theory, 2012, 58(11): 6942--6958.

\bibitem{Polak13}
M. Polak, U. Romanczuk, V. Ustimenko, A. Wroblewska,
On the applications of extremal graph theory to coding theory and cryptography,
Electron. Notes Discrete Math. 43 (2013) 329--342.

\bibitem{Cheng14}
X. Cheng, W. Chen, Y. Tang,
On the girth of the bipartite graph $D(k,q)$,
Discrete Math. 335 (2014) 25--34.

\bibitem{Cheng16}
X. Cheng, W. Chen, Y. Tang,
On the conjecture for the girth of the bipartite graph $D(k,q)$,
Discrete Math. 339 (2016) 2384--2392.

\bibitem{Sun17}S. Sun, On some families of algebraically defined graphs, Ph.D. Thesis, University of Delaware, 2017, http://udel.edu/shuying/girthDkq.pdf.

\bibitem{LazebnikSunWang17}
F. Lazebnik, S. Sun, Y. Wang,
Some families of graphs, hypergraphs and digraphs defined by systems of equations: a survey,
Lect. Notes Semin. Interdiscip. Mat. 14 (2017) 105--142.

\bibitem{Dehghan20}
A. Dehghan, A. H. Banihashemi,
Counting Short Cycles in Bipartite Graphs: A Fast Technique/Algorithm and a Hardness Result,
IEEE Trans. Commun. 68(3) (2020) 1378--1390.


\bibitem{XCT2023}
M. Xu, X. Cheng, Y. Tang,
Girth of the algebraic bipartite graph $D(k,q)$,
Finite Fields Appl. 92 (2023) 102296.

\bibitem{XCT23}
M. Xu, X. Cheng, Y. Tang,
On the girth cycles of the bipartite graph $D(k,q)$,




























Discrete Math. 346 (2023) 113500.

\bibitem{CZTX26}X. Cheng, C. Zhang, Y. Tang,
The backtrackless circuits of $D(k,q)$ obtained by lifting, submitted to Discrete Math.
\end{thebibliography}
\end{document}